\catcode`\@=11

\font\tenmsa=msam10 \font\sevenmsa=msam7 \font\fivemsa=msam5
\font\tenmsb=msbm10
\font\sevenmsb=msbm7 \font\fivemsb=msbm5 \newfam\msafam \newfam\msbfam
\textfont\msafam=\tenmsa \scriptfont\msafam=\sevenmsa
\scriptscriptfont\msafam=\fivemsa \textfont\msbfam=\tenmsb
\scriptfont\msbfam=\sevenmsb \scriptscriptfont\msbfam=\fivemsb

\def\hexnumber@#1{\ifnum#1<10 \number#1\else \ifnum#1=10 A\else\ifnum#1=11
 B\else\ifnum#1=12 C\else \ifnum#1=13 D\else\ifnum#1=14 E\else\ifnum#1=15
 F\fi\fi\fi\fi\fi\fi\fi}

\def\msa@{\hexnumber@\msafam} \def\msb@{\hexnumber@\msbfam}
\mathchardef\boxdot="2\msa@00 \mathchardef\boxplus="2\msa@01
\mathchardef\boxtimes="2\msa@02 \mathchardef\square="0\msa@03
\mathchardef\blacksquare="0\msa@04 \mathchardef\centerdot="2\msa@05
\mathchardef\lozenge="0\msa@06 \mathchardef\blacklozenge="0\msa@07
\mathchardef\circlearrowright="3\msa@08 \mathchardef\circlearrowleft="3\msa@09
\mathchardef\rightleftharpoons="3\msa@0A
\mathchardef\leftrightharpoons="3\msa@0B \mathchardef\boxminus="2\msa@0C
\mathchardef\Vdash="3\msa@0D \mathchardef\Vvdash="3\msa@0E
\mathchardef\vDash="3\msa@0F \mathchardef\twoheadrightarrow="3\msa@10
\mathchardef\twoheadleftarrow="3\msa@11 \mathchardef\leftleftarrows="3\msa@12
\mathchardef\rightrightarrows="3\msa@13 \mathchardef\upuparrows="3\msa@14
\mathchardef\downdownarrows="3\msa@15 \mathchardef\upharpoonright="3\msa@16
 \mathchardef\downharpoonright="3\msa@17
\mathchardef\upharpoonleft="3\msa@18 \mathchardef\downharpoonleft="3\msa@19
\mathchardef\rightarrowtail="3\msa@1A \mathchardef\leftarrowtail="3\msa@1B
\mathchardef\leftrightarrows="3\msa@1C \mathchardef\rightleftarrows="3\msa@1D
\mathchardef\Lsh="3\msa@1E \mathchardef\Rsh="3\msa@1F
\mathchardef\rightsquigarrow="3\msa@20
\mathchardef\leftrightsquigarrow="3\msa@21 \mathchardef\looparrowleft="3\msa@22
\mathchardef\looparrowright="3\msa@23 \mathchardef\circeq="3\msa@24
\mathchardef\succsim="3\msa@25 \mathchardef\gtrsim="3\msa@26
\mathchardef\gtrapprox="3\msa@27 \mathchardef\multimap="3\msa@28
\mathchardef\therefore="3\msa@29 \mathchardef\because="3\msa@2A
\mathchardef\doteqdot="3\msa@2B 
\mathchardef\traceiangleq="3\msa@2C \mathchardef\precsim="3\msa@2D
\mathchardef\lesssim="3\msa@2E \mathchardef\lessapprox="3\msa@2F
\mathchardef\eqslantless="3\msa@30 \mathchardef\eqslantgtr="3\msa@31
\mathchardef\curlyeqprec="3\msa@32 \mathchardef\curlyeqsucc="3\msa@33
\mathchardef\preccurlyeq="3\msa@34 \mathchardef\leqq="3\msa@35
\mathchardef\leqslant="3\msa@36 \mathchardef\lessgtr="3\msa@37
\mathchardef\backprime="0\msa@38 \mathchardef\risingdotseq="3\msa@3A
\mathchardef\fallingdotseq="3\msa@3B \mathchardef\succcurlyeq="3\msa@3C
\mathchardef\geqq="3\msa@3D \mathchardef\geqslant="3\msa@3E
\mathchardef\gtrless="3\msa@3F \mathchardef\sqsubset="3\msa@40
\mathchardef\sqsupset="3\msa@41
\mathchardef\trianglelefteq="3\msa@45 \mathchardef\bigstar="0\msa@46
\mathchardef\between="3\msa@47 \mathchardef\blacktriangledown="0\msa@48
\mathchardef\blacktriangleright="3\msa@49
\mathchardef\blacktriangleleft="3\msa@4A
\mathchardef\blacktriangle="0\msa@4E \mathchardef\triangledown="0\msa@4F
\mathchardef\eqcirc="3\msa@50 \mathchardef\lesseqgtr="3\msa@51
\mathchardef\gtreqless="3\msa@52 \mathchardef\lesseqqgtr="3\msa@53
\mathchardef\gtreqqless="3\msa@54 \mathchardef\Rrightarrow="3\msa@56
\mathchardef\Lleftarrow="3\msa@57 \mathchardef\veebar="2\msa@59
\mathchardef\barwedge="2\msa@5A \mathchardef\doublebarwedge="2\msa@5B
\mathchardef\angle="0\msa@5C \mathchardef\measuredangle="0\msa@5D
\mathchardef\sphericalangle="0\msa@5E \mathchardef\varpropto="3\msa@5F
\mathchardef\smallsmile="3\msa@60 \mathchardef\smallfrown="3\msa@61
\mathchardef\Subset="3\msa@62 \mathchardef\Supset="3\msa@63
\mathchardef\Cup="2\msa@64  \mathchardef\Cap="2\msa@65
 \mathchardef\curlywedge="2\msa@66
\mathchardef\curlyvee="2\msa@67 \mathchardef\leftthreetimes="2\msa@68
\mathchardef\rightthreetimes="2\msa@69 \mathchardef\subseteqq="3\msa@6A
\mathchardef\supseteqq="3\msa@6B \mathchardef\bumpeq="3\msa@6C
\mathchardef\Bumpeq="3\msa@6D \mathchardef\lll="3\msa@6E 
\mathchardef\ggg="3\msa@6F  \mathchardef\circledS="0\msa@73
\mathchardef\pitchfork="3\msa@74 \mathchardef\dotplus="2\msa@75
\mathchardef\backsim="3\msa@76 \mathchardef\backsimeq="3\msa@77
\mathchardef\complement="0\msa@7B \mathchardef\intercal="2\msa@7C
\mathchardef\circledcirc="2\msa@7D \mathchardef\circledast="2\msa@7E
\mathchardef\circleddash="2\msa@7F \def\ulcorner{\delimiter"4\msa@70\msa@70 }
\def\urcorner{\delimiter"5\msa@71\msa@71 }
\def\llcorner{\delimiter"4\msa@78\msa@78 }
\def\lrcorner{\delimiter"5\msa@79\msa@79 } \def\yen{\mathhexbox\msa@55 }
\def\checkmark{\mathhexbox\msa@58 } \def\circledR{\mathhexbox\msa@72 }
\def\maltese{\mathhexbox\msa@7A } \mathchardef\lvertneqq="3\msb@00
\mathchardef\gvertneqq="3\msb@01 \mathchardef\nleq="3\msb@02
\mathchardef\ngeq="3\msb@03 \mathchardef\nless="3\msb@04
\mathchardef\ngtr="3\msb@05 \mathchardef\nprec="3\msb@06
\mathchardef\nsucc="3\msb@07 \mathchardef\lneqq="3\msb@08
\mathchardef\gneqq="3\msb@09 \mathchardef\nleqslant="3\msb@0A
\mathchardef\ngeqslant="3\msb@0B \mathchardef\lneq="3\msb@0C
\mathchardef\gneq="3\msb@0D \mathchardef\npreceq="3\msb@0E
\mathchardef\nsucceq="3\msb@0F \mathchardef\precnsim="3\msb@10
\mathchardef\succnsim="3\msb@11 \mathchardef\lnsim="3\msb@12
\mathchardef\gnsim="3\msb@13 \mathchardef\nleqq="3\msb@14
\mathchardef\ngeqq="3\msb@15 \mathchardef\precneqq="3\msb@16
\mathchardef\succneqq="3\msb@17 \mathchardef\precnapprox="3\msb@18
\mathchardef\succnapprox="3\msb@19 \mathchardef\lnapprox="3\msb@1A
\mathchardef\gnapprox="3\msb@1B \mathchardef\nsim="3\msb@1C
\mathchardef\napprox="3\msb@1D
\mathchardef\nsupseteqq="3\msb@23 \mathchardef\subsetneqq="3\msb@24
\mathchardef\supsetneqq="3\msb@25
\mathchardef\supsetneq="3\msb@29 \mathchardef\nsubseteq="3\msb@2A
\mathchardef\nsupseteq="3\msb@2B \mathchardef\nparallel="3\msb@2C
\mathchardef\nmid="3\msb@2D \mathchardef\nshortmid="3\msb@2E
\mathchardef\nshortparallel="3\msb@2F \mathchardef\nvdash="3\msb@30
\mathchardef\nVdash="3\msb@31 \mathchardef\nvDash="3\msb@32
\mathchardef\nVDash="3\msb@33 \mathchardef\ntrianglerighteq="3\msb@34
\mathchardef\ntrianglelefteq="3\msb@35 \mathchardef\ntriangleleft="3\msb@36
\mathchardef\ntriangleright="3\msb@37 \mathchardef\nleftarrow="3\msb@38
\mathchardef\nrightarrow="3\msb@39 \mathchardef\nLeftarrow="3\msb@3A
\mathchardef\nRightarrow="3\msb@3B \mathchardef\nLeftrightarrow="3\msb@3C
\mathchardef\nleftrightarrow="3\msb@3D \mathchardef\divideontimes="2\msb@3E
\mathchardef\varnothing="0\msb@3F \mathchardef\nexists="0\msb@40
\mathchardef\mho="0\msb@66 \mathchardef\thorn="0\msb@67
\mathchardef\beth="0\msb@69 \mathchardef\gimel="0\msb@6A
\mathchardef\daleth="0\msb@6B \mathchardef\lessdot="3\msb@6C
\mathchardef\gtrdot="3\msb@6D \mathchardef\ltimes="2\msb@6E
\mathchardef\rtimes="2\msb@6F \mathchardef\shortmid="3\msb@70
\mathchardef\shortparallel="3\msb@71 \mathchardef\smallsetminus="2\msb@72
\mathchardef\thicksim="3\msb@73 \mathchardef\thickapprox="3\msb@74
\mathchardef\approxeq="3\msb@75 \mathchardef\succapprox="3\msb@76
\mathchardef\precapprox="3\msb@77 \mathchardef\curvearrowleft="3\msb@78
\mathchardef\curvearrowright="3\msb@79 \mathchardef\digamma="0\msb@7A
\mathchardef\varkappa="0\msb@7B \mathchardef\hslash="0\msb@7D
\mathchardef\hbar="0\msb@7E \mathchardef\backepsilon="3\msb@7F
\def\Bbb{\ifmmode\let\next\Bbb@\else
\def\next{\errmessage{Use \string\Bbb\space only in math mode}}\fi\next}
\def\Bbb@#1{{\Bbb@@{#1}}} \def\Bbb@@#1{\fam\msbfam#1}

\catcode`\@=\active

\def\del{\partial}

 \def\CC{\hbox{{$\cal C$}}}
 
\def\CL{\hbox{{$\cal L$}}}

\def\CR{\hbox{{$\cal R$}}}

\def\CQ{\hbox{{$\cal Q$}}}

\def\cg{\hbox{{\sl g}}} 

\def\lform{\hbox{$\sqcup$}\llap{\hbox{$\sqcap$}}}

\def\h{{{1\over2}}}

\def\R{{\Bbb R}}
\def\C{{\Bbb C}}

\def\eps{{\epsilon}}

\def\dcross{{\bowtie}}
\def\codcross{{\blacktriangleright\!\!\blacktriangleleft}}
\def\rbiprod{{\cdot\kern-.33em\triangleright\!\!\!<}}
\def\lbiprod{{>\!\!\!\triangleleft\kern-.33em\cdot\, }}

\def\tens{\mathop{\otimes}}

\def\la{{\triangleright}}\def\ra{{\triangleleft}}

\def\extd{{{\rm d}}}

\def\isom{{\cong}}

\def\Ad{{\rm Ad}}

\def\id{{\rm id}}

\def\Lin{{\rm Lin}}

\def\<{\langle}
\def\>{\rangle}

\def\equad{\kern -1.7em}
\def\qqquad{\qquad\quad}

\def\eqn#1#2{\begin{equation}#2\label{#1}\end{equation}}

\def\haj#1{{\mathaccent20 {#1}}}

\def\o{{}_{\scriptscriptstyle(1)}}
\def\t{{}_{\scriptscriptstyle(2)}}
\def\th{{}_{\scriptscriptstyle(3)}}

\def\bo{{}^{\bar{\scriptscriptstyle(1)}}}
\def\bt{{}^{\bar{\scriptscriptstyle(2)}}}

\def\Ruo#1{{\CR^{\scriptscriptstyle(1)}_{#1}}}
\def\Rut#1{{\CR^{\scriptscriptstyle(2)}_{#1}}}

\def\uo{{{}^{\scriptscriptstyle(1)}}}
\def\ut{{{}^{\scriptscriptstyle(2)}}}

\def\umo{{{}^{\scriptscriptstyle-(1)}}}
\def\umt{{{}^{\scriptscriptstyle-(2)}}}

\def\text#1{\mbox{\rm #1}}
\def\note#1{}

\def\blacksquare{{\lform}}
\def\frac#1#2{{{#1\over#2}}}

\def\proof{\goodbreak\noindent{\bf Proof\quad}}

\def\endproof{{\ $\lform$}\bigskip }

\def\align#1{\begin{eqnarray*}#1\end{eqnarray*}}

\def\ceqn#1#2{\begin{equation}\label{#1}
\begin{array}{c}#2\end{array}\end{equation}}

\def\vecu{{\bf u}}

\documentstyle[11pt]{article}
\newtheorem{lemma}{Lemma}[section] \newtheorem{propos}[lemma]{Proposition}
 \newtheorem{theorem}[lemma]{Theorem}

\begin{document}\baselineskip 21pt

{\ }\qquad\qquad \hskip 4.3in 
\vspace{-.2in}

\begin{center} {\LARGE CLASSIFICATION OF BICOVARIANT}\\ {\LARGE DIFFERENTIAL
CALCULI}
\\ \baselineskip 13pt{\ }
{\ }\\
S. Majid\footnote{Royal Society University Research Fellow and Fellow of
Pembroke College, Cambridge}\\
{\ }\\
Department of Mathematics, Harvard University\\
Science Center, Cambridge MA 02138, USA\footnote{During the calendar years 1995
+ 1996}\\
+\\
Department of Applied Mathematics \& Theoretical Physics\\
University of Cambridge, Cambridge CB3 9EW\\
\end{center}

\begin{center}
July -- August 1996
\end{center}
\vspace{10pt}
\begin{quote}\baselineskip 13pt
\noindent{\bf Abstract}
We show that the bicovariant first order differential calculi on a
factorisable semisimple quantum group are in 1-1 correspondence with
irreducible representations $V$ of the quantum group enveloping
algebra. The corresponding calculus is constructed and has dimension
${\rm dim} V^2$. The differential calculi on a finite group algebra $\C
G$ are also classified and shown to be in correspondence with pairs
consisting of an irreducible representation $V$ and a continuous
parameter in $\C P^{\dim V -1}$. They have dimension $\rm dim V$. For a
classical Lie group we obtain an infinite family of non-standard
calculi. General constructions for bicovariant calculi and their
quantum tangent spaces are also obtained.


\end{quote}
\baselineskip 21pt

\section{Introduction}

One of the first steps in non-commutative geometry of the kind coming out of
quantum
groups is the choice of `first order differential calculus' or `cotangent
bundle'. Only once this
 is fixed can one   begin to do gauge theory  \cite{BrzMa:gau} or make other
geometrical
constructions.  When the quantum space in question is a  quantum group, $A$, it
is natural to require that the differential calculus is covariant under left
and right translations. Thus, we require
$\Gamma\equiv\Omega^1(A)$, $\extd:A\to \Gamma$  such that

1. $\Gamma$ is an $A$-bimodule

2. $\Gamma$ is an $A$-cobimodule, with coactions $\Delta_L:\Gamma\to A\tens
\Gamma$
and $\Delta_R:\Gamma\to \Gamma\tens A$ bimodule maps

3. $\extd:A\to \Gamma$ is a bicomodule map

4. $\extd(ab)=(\extd a)b+a(\extd b)$ for all $a,b\in A$

5. $\Gamma={\rm span}\{a\extd b|\ a,b\in A\}$

Here, a bicomodule is like a bimodule but with arrows reversed, i.e. a pair of
commuting coactions
$\Delta_L,\Delta_R$. A morphism of differential calculi means a bimodule and
bicomodule map forming a commutative triangle with the $\extd$ maps. These are
the natural axioms proposed some years ago by Woronowicz\cite{Wor:dif}. The
axiom 5. here is not essential;   it specifies that the calculus is
irreducible, and is assumed throughout the paper. By now, several
examples are known, and there is also case-by-case classification for several
families of quantum groups in \cite{SchSch:bic}. The class of `inner'
bicovariant
calculi has also been introduced \cite{BrzMa:bic}. The complete classification
of the
possible calculi on a general quantum group  remains, however, open.

In Section~4 of the present paper we present a complete solution to this
classification
problem under the strict assumption of a semisimple factorisable quantum group.
The standard $q$-deformed function algebras $G_q$ are essentially factorisable
in the sense that they are factorisable up to suitable localisations or when
working
over formal power-series. In this case our algebraic result (a) constructs a
calculus
of on $G_q$ of dimension $({\rm dim} V)^2$ for each irreducible representation
$V$ of $U_q(\cg)$
and (b) indicates that these are the only `generic' possibilities in the sense
of
extending to localisations or to working over formal power-series in the
deformation
parameter.

We begin in Section~2 with a clarification of the role of the quantum double in
classifying calculi. This is well-known or implicit from \cite{Wor:dif} but
appears
to remain of current interest; see \cite{BGMST:dou}. In fact, calculi
correspond to subrepresentations of a given quantum-double module, a
result which is somewhat different from that recently presented in
\cite{BGMST:dou} (these authors did not impose the optional
irreducibility axiom~5 above, and hence have a more complicated
result). We reformulate the theory with these quantum double
subrepresentations or `quantum tangent spaces' as the starting point
and then develop some preliminary general results from this point of
view.

In Section~3 we apply these results to
the complete classification for $A=\C G$ the group algebra over a finite group,
as well as recovering the known classification for the function algebra
$A=\C(G)$. We also comment on the case where $G$ is a Lie group.

Section~4 presents the main result of the paper; Section~5 concludes with a
brief discussion of the open problem for the classification of higher order
differential calculi (or exterior algebras). We also mention the braided
version of the results in the present paper, to be presented in detail
elsewhere\cite{Ma:bdif}. We work over $\C$ for convenience, but all abstract
Hopf-algebraic results work over any ground field or, with care, over a ring
such as $\C[[\hbar]]$.

\subsection*{Acknowledgments} Most of work for this paper and some of the
writing up was done at the meeting on
Representation Theory and Mathematical Physics at the Erwin Schroedinger
Institute during my visit
May-June 1996. I would like thank the organisers I. Penkov, J. Wolf and the
director P. Michor
for their support.

\section{Role of the quantum double}

In this preliminary section we clarify the role of the quantum double in
classifying bicovariant differential
calculi. It is needed for our main result in  Section~4. This role of the
quantum double should be known to experts, but we have not found a suitably
explicit treatment elsewhere. Moreover,
this use of the quantum double is different from the one recently presented in
\cite{BGMST:dou}, hence it would appear necessary to emphasise it here and
explain the relation with that work.  The use of a braiding to describe the
derivation property of partial derivatives in Proposition~2.3,   the emphasis
on quantum double subrepresentations, and the resulting applications  at the
end of the section (such as the `mirror' operation) appear to be novel aspects
of our formulation.

Let $H$ be a Hopf algebra non-degenerately dually paired with $A$. The Drinfeld
quantum double\cite{Dri} is the double cross product Hopf algebra $H\dcross
A^{\rm op}$ built on
$H\tens A$ with the product
\eqn{double}{ (h\tens a)(g\tens b)=hg\t\tens ba\t
\<g\o,a\o\>\<g\th,Sa\th\>,\quad a,b\in A,\ h,g\in H }
and tensor product unit and coalgebra. This is the formulation from
\cite{Ma:dou}. We use here (and throughout) the notations and
conventions   from \cite{Ma:book}. Thus, $\Delta h=h\o\tens h\t$ is the
coproduct, $S$ is the antipode (which we assume for convenience to be
invertible), and $\<\ , \ \>$ is the pairing between $H$ and $A$. We denote the
counit
of any of our Hopf algebras by $\eps$.

The quantum double has a formal quasitriangular structure $\CR=\sum_a f^a\tens
e_a$ where $\{e_a\}$ is a basis of $H$ and $\{f^a\}$ a dual basis of $A$.
Although formal, this does lead to a braiding
$\Psi$ among suitable representations. To make this precise, we define a
representation of the quantum double of $H$ to be $H$-regular if the action of
$A\subset H\dcross A^{\rm op}$ is given by evaluation against
a (left) coaction of $H$. It is $A$-regular if the action of $H$ is given by
evaluation against a (right) coaction of $A$. If $V$ is $A$-regular or $W$ is
$H$-regular then $\Psi:V\tens W\to W\tens V$ is a well-defined operator. Thus,
$\Psi(v\tens w)=\sum_a e_a\la w\tens f^a\la v=\sum_a e_a\la w\tens
\<f^a,v\bo\>v\bt=v\bo\la w\tens v\bt$ in the first case, where $v\mapsto
v\bo\tens v\bt$ (with summation understood) is the assumed coaction $V\to
H\tens V$. Similarly in the second case.

Woronowicz in \cite{Wor:dif} observed that first order bicovariant calculi are
in
1-1 correspondence with $\Ad$-invariant `ideals' in $\ker\eps$. More precisely,
they correspond to   quotients of $\ker\eps\subset A$ by subspaces $M$ which
are stable under the  action and
coaction
\eqn{quotcross}{ a\la v=av,\quad \Ad(v)= v\o Sv\th\tens  v\t}
on $v\in\ker\eps$. Equivalently, they correspond to quotients $V$ to which this
action and coaction descend. The corresponding calculus is $\Gamma=V\tens A$
with tensor product action and coaction from the left and
trivial action and coaction on $V$ from the right (here we take the left and
right actions and coactions on $A$
defined by its product and coproduct). In addition, $\extd a= a\o\tens
a\t-1\tens a$, where $a\o$ is projected to $V$. This is the most general form
for a bicovariant calculus up to isomorphism. The case $V=\ker\eps\subset A$ is
called the `universal' first order calculus.

As a first step, we can write all coactions of $A$ as actions of $H$. Then
(\ref{quotcross}) becomes
\eqn{quotrep}{  a\la v=av,\quad h\la v=\<h,v\o Sv\th\> v\t}
and calculi correspond to quotients of which are equivariant under these
actions. The quantum double is not needed to classify calculi here, but in fact
these two actions do fit together
to form a representation of $A\dcross H^{\rm op}$, the quantum double of $A$.
This fact allows one to deduce,
for example, the canonical braiding $\sigma$ in \cite{Wor:dif} from the
quantum double quasitriangular structure, which would otherwise have to
be introduced by hand. This was explained in \cite{Ma:sol}. In this context, it
is natural also to reformulate the bicomodules $\Delta_L,\Delta_R$ in the
axioms 2. and 3. of a bicovariant calculus  as an $H^{\rm op}$-bimodule by
evaluation against the coactions.
In principle, requiring an $H^{\rm op}$-bimodule could be slightly more general
when $H$ is infinite-dimensional.

\begin{lemma} The quantum double $H\dcross A^{\rm op}$ acts on $\ker\eps\subset
H$ by
\[ h\la x=h\o x Sh\t,\quad a\la x=\<a,x\o\>x\t-\<a,x\>1\]
\end{lemma}
\proof The quantum double has a well-known `Schroedinger' representation on
$H$\cite{Ma:book}
by the quantum adjoint action and by the coregular `differentiation'
representation. This induces the
action stated on  $\ker\eps$ via the projection $\Pi(h)=h-1\eps(h)$ as a
morphism $H\to \ker\eps$, i.e. it is
easy to see that it is indeed an action of the quantum double on $\ker\eps$ and
that $\Pi$ is
an intertwiner. Also,  we can identify the linear space $\ker\eps\subset A$
with $A/\C$ (the
quotient by the 1-dimensional vector space spanned by the unit element); for
any element in $A/\C$ there is a unique representative in $\ker\eps\subset A$.
In terms of $A/\C$ the action in (\ref{quotrep}) is
$a\la v=av-\eps(v)a$ and  $h\la v=\<v\o Sv\th\> v\t$. The action stated in the
lemma is the natural right action of $A\dcross H^{\rm op}$ on $\ker\eps\subset
H$ dual to this action on $A/\C$, viewed as a left action of  the quantum
double of $H$. \endproof

We are now ready to make a further reformulation which, when the bicovariant
calculus is
finite-dimensional as an $A$-module, is strictly equivalent by dualizing $V$.

\begin{propos} Finite-dimensional bicovariant calculi are in 1-1 correspondence
with
subrepresentations $L\subseteq \ker\eps\subset H$ of the quantum double
representation in Lemma~2.1.
\[ \Gamma=\Lin(L,A),\qquad (\extd a)(x)=\<x,a\o\>a\t \]
\[ (a\cdot\gamma)(x)=a\t\gamma(a\o\la x),\quad (\gamma\cdot a)(x)=\gamma(x)a\]
\[ (h\cdot\gamma)(x)=\<h\t,\gamma(h\o\la x)\o\>\gamma(h\o\la x)\t,\quad
(\gamma\cdot h)(x)=\gamma(x)\o \<\gamma(x)\t,h\>\]
for all $\gamma\in\Lin(L,A)$. The vector space $L$ is called  the {\em quantum
tangent space} associated to the
calculus.
\end{propos}
\proof A quotient of $A/\C$ which is equivariant under the quantum double
action corresponds under dualisation to
a subspace of $\ker\eps\subset H$ which is stable under the action in
Lemma~2.1.  Thus, this is a dual formulation of
the correspondence (\ref{quotrep}). Also, $\extd a=a\o\tens a\t-1\tens
a=a\o\tens a\t$ when we work with $V$ as a quotient of $A/\C$, which leads to
the form shown for $\extd$. It is easy to verify directly that the structures
shown indeed provide
a first order bicovariant differential calculus given $L$. Conversely, in the
finite-dimensional case, we define $L=V^*$ where $V$ is the invariant part of
$\Gamma$ under the usual correspondence in (\ref{quotrep}).  In terms of the
ideal $M$ which defines the bicovariant calculus via (\ref{quotcross}), the
corresponding quantum double subrepresentation is $L=\{x\in\ker\eps|\
\<x,a\>=0\ \forall a\in M\}$. \endproof

The correspondence in the proposition is contragradient i.e. morphisms of
calculi $\Gamma_1\to\Gamma_2$ correspond to inclusions $L_2\hookrightarrow L_1$
of quantum double subrepresentations.  Only inclusions are allowed here,
corresponding to all morphisms of calculi being of the form $\Gamma_2$ a
quotient of $\Gamma_1$. In the infinite-dimensional case every bicovariant
calculus continues to define a subrepresentation $L$, and, conversely, a
subrepresentation $L$ continues to provide a bicovariant first order calculus
in our slightly generalised sense where an action of $H^{\rm op}$ replaces the
coactions $\Delta_L,\Delta_R$. It is this final reformulation in terms quantum
tangent spaces which we will use; by definition a quantum tangent space $L$ is
a subrepresentation of $\ker\eps\subset H$ under the action of the quantum
double of $H$, and it is these which we will actually classify in the present
paper. Indeed, quantum tangent spaces have many nice properties, making them an
equally good starting point for differential calculus.

\begin{propos} For each $x\in L$, we define the `braided derivation'
\[ \del_x:A\to A,\qquad \del_x(a)=(\extd a)(x)\]
This obeys
\[ \del_x(ab)=(\del_xa)b+  \Psi(a\tens x)\umt\del_{\Psi(a\tens x)\umo}b\]
where $\Psi:L\tens A\to A\tens L$ is the quantum double braiding between $A,L$
as quantum double modules, with inverse $\Psi^{-1}(a\tens x)$ denoted
explicitly by $\Psi(a\tens x)\umo\tens\Psi(a\tens x)\umt$.
\end{propos}
\proof  We start with the identity
\[\del_x(ab) =\<x,(ab)\o\>(ab)\t= \<x\o,a\o\>\<x\t,b\o\>a\t b\t\qqquad\qqquad\]
\[ =\<x,a\o\> a\t b+ \<a\o\la x,b\o\> a\t b\t =\del_x(a)b+a\t\del_{a\o\la
x}(b)\]
based on the definition of $\del_x$ and the action in Lemma~2.1. On the other
hand, $A$ is
a quantum double module by
\eqn{Arep}{h\la a=\<Sh,a\o\>a\t,\quad b\la a=(S^{-1}b\t)a b\o,}
which is the conjugate (dual) of the Schroedinger representation of the quantum
double on $H$. It is $H$-regular, so that the braiding $\Psi$ is well-defined.
We compute it easily as
\eqn{Psixa}{\Psi(x\tens a)= e_a\la a\tens f^a\la x= a\t\tens Sa\o\la x}
with inverse $\Psi^{-1}(a\tens x)=a\o\la x\tens a\t$, which we put into the
above identity.
\endproof

Because the subspace $L$ is stable under the quantum adjoint action, it is
tempting
to restrict the latter to $L$ as a `quantum Lie bracket' $[\ ,\ ]=\Ad:L\tens
L\to L$.
The use of $\Ad$ as `quantum Lie bracket' has been discussed in \cite{Wor:dif}
from this point of view and independently from another point of view in
\cite{Ma:skl}  (where the `quantum Lie bracket structure constants' for the
$l^+Sl^-$ generators
of $U_q(\cg)$ were computed in R-matrix form). The only content here
comes from the
identities
\eqn{Hjac}{[x,[y,z]]=[[x\o,y],[x\t,z]]{},\quad [x\o,y]x\t=xy}
which hold in any quantum group when $[\ ,\ ]=\Ad$ is the left adjoint action.

\begin{propos} The `quantum Lie bracket' on $L$ defined by $\Ad$ obeys
\[ [x,[y,z]]=[[x,y],z]+[\ ,[\ ,z]]\circ\Psi(x\tens y),\quad [x,y]=xy-\cdot
\Psi(x\tens y) \]
where $\Psi$ is the quantum double braiding between $L,L$ as quantum double
modules.
\end{propos}
\proof We again use the formula $\Psi(x\tens y)=e_a\la y\tens f^a\la x$;  our
action of the quantum double of $H$ is regular and we still have a well-defined
operator
\eqn{Psixy}{\Psi(x\tens y)=\Ad_{e_a}(y)\tens (\<f^a,x\o\>x\t - \<f^a,x\>1)=
[x\o,y]\tens x\t - [x,y]\tens 1.}
Then (\ref{Hjac}) can be trivially rewritten in the form stated by
eliminating the coproduct in favour of $\Psi$ in these equations. \endproof

There are, however, some fundamental problems to be overcome before one could
call this
vector space $L$ with $[\ ,\ ]$ some kind of `quantum Lie algebra'. These
fundamental
problems have been explained in \cite{Ma:skl} and force one to the braided
version\cite{Ma:lie} where these problems are resolved:

1. Although  `enveloping algebra-like'  relations $[x,y]=xy-\cdot \Psi(x\tens
y)$ hold in $H$, we do not
know that $L$ generates $H$. Even if it does, we do not know that these are the
only relations in $H$. Indeed, for $U_q(\cg)$ they are not. So $H\ne U(L)$ as
generated by $L$ and such relations.

2. Even if $H$ were to be generated in some way from $L$,
we are not able to recover the coproduct of $H$ in this way. Indeed, for
$U_q(\cg)$ the coproduct of $H$ does not have any simple form on
$L$ and hence cannot be generated in some canonical way. Without this, $U(L)$
is only an
algebra and not a Hopf algebra or bialgebra. Equivalently,  one cannot tensor
product representations of $L$ in any natural way, which makes it useless as a
`Lie algebra'.

In the case where $H$ is quasitriangular, there is a `transmutation
theory'\cite{Ma:bra} which
converts $H$ to a braided group. It also converts $L$ to a `braided-Lie
algebra'
$\CL$. The linear maps $[\ ,\ ]$ are the same (the
braided adjoint action coincides with the quantum one) but the coalgebras
are different. For the standard calculus on $G_q$, the braided coproduct takes
a standard
matrix form on $\CL$ and there is a corresponding $U(\CL)$ as
a braided group (bialgebra in a braided category) generated from $\CL$. Thus,
the problems 1.-2. are
resolved at the price of working with the braided version of the theory.
For a general Hopf algebra $H$ and general $L$, however,  we do not really have
a `quantum Lie algebra' or braided-Lie
algebra at all. Hence we  prefer the term `quantum tangent space'   for the
subspace $L$.

Finally, we discuss the well-known correspondence between bicovariant
bimodules (i.e. simultaneous modules and comodules as in the above axioms 1.
and 2. for a calculus) and quantum double modules. Bicovariant bimodules are
known in the mathematics literature\cite{Nic:bia} as bi-Hopf modules, and
indeed correspond (in a standard way) to crossed modules\cite{Lyu:pri}. The
latter have been identified with quantum double
modules by the author in \cite{Ma:dou}. In view of such a correspondence, it
has
been argued in \cite{BGMST:dou} that bicovariant calculi are therefore in
1-1 correspondence with pairs $(V,\psi)$ where $V$ is a quantum double module
and
$\psi$ is an $\Ad$-equivariant one-cocycle on $A$ with values in $V$. The
corresponding
calculus is
\[ \Gamma=V\tens A,\quad \extd a=\psi(a\o)\tens a\t\]
It is easy to see that the Leibniz rule for $\extd$ indeed corresponds to the
cocycle condition
\[ \psi(ab)=a\la\psi(b)+\psi(a)\eps(b).\]
By contrast, the classification in Proposition~2.2 tells us that this point of
view, although interesting, is not so useful; not every quantum double module
$V$ is needed but only the subrepresentations  of {\em one particular quantum
double representation}.
This is because the above axiom~5. of a bicovariant calculus corresponds to
$\psi$
surjective. This condition was omitted in the analysis in \cite{BGMST:dou} and,
as soon as it is imposed, we see that the cocycle condition forces the action
of $A$
since every element of $V$ can be written as $\psi(b)$ for some $b$. Likewise,
surjectivity
of $\psi$ and its equivariance forces the action of $H$ on $V$ to be the image
under $\psi$ of the coadjoint one. Therefore, we are forced into the setting
above,
with $V=L^*$ and $\psi$ the quotienting map in (\ref{quotrep}) or the adjoint
of the
inclusion $L\subseteq\ker\eps$ in Proposition~2.2. This is why not all quantum
double modules are allowed. Moreover, it is not necessary to solve any cocycle
conditions explicitly; these take care of themselves in the specification of
the
submodule inclusion into $\ker\eps$.

We conclude the section with some general constructions for bicovariant calculi
and their quantum tangent spaces. Firstly, we recall that any element
$\alpha\in A$ which is invariant under the adjoint coaction $\Ad$ in
(\ref{quotcross}) can be used
to generate an ideal $A(\alpha-\eps(\alpha))$ to quotient $\ker\eps\subset A$
by. This class of bicovariant calculi can be called {\em inner type-I} because
the exterior derivative obeys
\eqn{inner}{ \eps(\alpha)\extd a=a\o\eps(\alpha)\tens a\t -1\eps(\alpha)\tens
a=a\o\alpha\tens a\t-\alpha\tens a=a\cdot(\alpha\tens 1)-(\alpha\tens 1)\cdot
a}
projected down to the quotient. The expression on the right is shown lifted up
to $A\tens A$ as a left $A$-module by the tensor product left action and a
right $A$-module by right multiplication in the second copy. A variant of this
construction was introduced in \cite{BrzMa:bic}, where we quotient by the
ideal $(\ker\eps).(\alpha
-(\eps(\alpha)+1))\subseteq\ker\eps\subset A$ and we have
\ceqn{innerII}{ \extd a=(a\o-\eps (a\o))\tens a\t=
a\o(\alpha-\eps(\alpha))\tens a\t -(\alpha-\eps(\alpha))\tens a
=a\cdot\omega(\alpha)-\omega(\alpha)\cdot a}
in the quotient, where $\omega(\alpha)=(\alpha-\eps(\alpha))\tens
1=(\extd \alpha\o)S\alpha\t \in \Gamma$. This class can be called {\em
inner type-II}. Since $\Ad$-invariant elements $\alpha$ are closed
under addition and multiplication, we have whole ring of bicovariant
differential calculi of either type. The standard calculi on $G_q$ were
already obtained in \cite{Jur:dif} as a quotient of an inner type-II
form (with $\alpha$ the $q$-trace), while \cite{BrzMa:bic} extended this
to a ring of calculi generated by elements $\alpha_1,\cdots,\alpha_r\in
G_q$ constructed through transmutation. Note that the cocycle point of
view  might suggest the more general notion of `coboundary'
differential calculus where $\psi$ is the Hocschild coboundary of
$\nu\in V$, i.e. $\psi(a)=a\la \nu-\eps(a)\nu$. However, when we again
add the surjectivity of $\psi$ we are forced to $\nu=\psi(\alpha)$ for
some $\alpha$ and we return to the class of inner type-II calculi from
\cite{BrzMa:bic}. This is simlar to (but different from) the discussion
in \cite{BGMST:dou}.

\begin{propos} The quantum tangent space for an inner type-I
bicovariant differential calculus defined by any non-trivial element
$\alpha\in A$ invariant under $\Ad$ in (\ref{quotcross}) is the quantum
double subrepresentation \[ L_\alpha=\{x\in \ker\eps\subset
H|\ x\o\<x\t,\alpha\>=x\eps(\alpha)\}.\] This has a canonical extension
$\tilde{L_\alpha}$ where the condition on $x$ is only required to hold
on evaluation against all $a\in\ker\eps\subset A$. Similarly, the quantum
tangent space for the inner type-II case is \[ L_{\alpha,1}= \{x\in
\ker\eps\subset H|\ \<x,a\alpha\>=\<x,a\>(\eps(\alpha)+1),\ \forall
a\in\ker\eps\subset A\}.\]
\end{propos}
\proof It is convenient to first identify $\ker\eps$ with $A/\C$ as in
the proof of Lemma~2.1. Then an inner type-I bicovariant  calculus has $V$ the
quotient by the image $A\la\alpha$ in $A/\C$.  Hence its dual consists
of the linear functionals   $x\in\ker\eps\subset H$ such that
$\<x,a\la\alpha\>=0$ for all $a$, i.e. such that
$\<x,a\alpha-a\eps(\alpha)\>=0$. This leads to the dual formulation; we
define $L_\alpha$ as stated and verify directly that it is stable under the
quantum double action in Lemma~2.1, which is a straightforward Hopf
algebra computation. The variant in which we require
$\<x,a\alpha\>=\<x,a\>\eps(\alpha)$ for all $a\in\ker\eps\subset A$ is
easily verified to also form a quantum double subrepresentation, and
defines $\tilde{L_\alpha}$. The type-II case is similar. \endproof

Note that we can define $L_{\alpha,\lambda}$ similarly, with
$\eps(\alpha)+\lambda$ in place of $\eps(\alpha)+1$ and then obtain
$\lambda \extd a=a\cdot\omega(\alpha)-\omega(\alpha)\cdot a$ as in
\cite{BrzMa:bic}.  All non-zero $\lambda$ are equivalent to the inner
type-II construction via $L_{\alpha,\lambda}=L_{\lambda^{-1}\alpha,1}$,
while $L_{\alpha,0}=\tilde{L_\alpha}$.  More generally, we can
restrict the condition $\forall a\in\ker\eps\subset A$ in
by requiring only $a\in M$, where $M\subseteq\ker\eps\subset A$ is the
quotienting ideal for any given
bicovariant calculus. This gives a 1-parameter family of new calculi with
quotienting ideal $M\cdot(\alpha-(\eps(\alpha)+\lambda))$, and is the general
idea behind the above constructions.

The representation-theoretic point of view suggests, however, a different
type of general construction for any Hopf algebra. Namely, pick any element
$x\in\ker\eps\subset H$ and define $L=(H\dcross A^{\rm op})\la x$, the image of
$x$ under the quantum double action. It evidently forms a subrepresentation of
$\ker\eps$ and hence by Proposition~2.2 it defines a bicovariant differential
calculus. More generally, the image of any left ideal of the quantum double
acting on any $x\in \ker\eps\subset H$ will be a subrepresentation. An
interesting special case of this idea is the following:

\begin{propos} Let $c\in H$ be any non-trivial central element. There is an
associated bicovariant differential calculus with
\[ L_c={\rm span}\{ x_a=\<a,c\o\> c\t - \<a,c\>1|\ a\in \ker\eps\subset A\}\]
\[\del_{x_a}(b)=\<ab\o,c\>b\t - \<a,c\>b,\quad\Psi^{-1}(b\tens x_a)=
x_{ab\o}\tens b\t\]
\[ [x_a,x_b]=x_{b\t}\<aSb\o)b\th,c\>-x_b\<a,c\>,\quad \Psi(x_a\tens
x_b)=x_{b\t}\tens x_{a(Sb\o)b\th}\]
for all $b\in A$ in the middle line. We say that the quantum tangent space
$L_c$ or its corresponding bicovariant differential calculus is {\em centrally
generated}. It has a canonical extension $\tilde{L_c}$ spanned by $\{x_a\}$ for
all $a\in A$.
\end{propos}
\proof The quantum double can also be written as $A^{\rm op}\dcross H$, i.e.
every element can be written uniquely in the form $\sum_i a_ih_i$ with $a_i\in
A$ and $h_i\in H$. A central element $c$ is precisely an element for which
$h\la c=\eps(h)c$ for all $h\in H$. Hence the image of $x=c-\eps(c)$ under the
quantum double action reduces to the image of the action of $A$ in Lemma~2.1.
Thus $\tilde{L_c}={\rm span}\{x_a=a\la (c-\eps(c))| a\in A\}$ is a
subrepresentation under the quantum double. We then restrict the allowed
$\{x_a\}$ to $a\in \ker\eps\subset A$. It is easy to check that this still
defines a quantum double subrepresentation, which is the one stated. It can
sometimes coincide with $\tilde{L_c}$.

The calculation of $\del_{x_a}$ and $\Psi^{-1}$ from Proposition~2.3 is
trivial. For the quantum Lie bracket in Proposition~2.4, we note first the
$\Ad$-invariance identity
\eqn{casinv}{ c\o\tens h\o c\t Sh\t =(Sh\o)c\o h\t\tens c\t,\quad \forall
h\in H}
which holds for any central element $c$. Then
\align{\Ad_h(x_a)\equad &&=\<a,c\o\>\Ad_h(c\t)-\eps(h)\<a,c\>=\<a,(Sh\o)c\o
h\t\>c\t-\eps(h)\<a,c\>\\
&&=x_{a\t}\<h,(Sa\o)a\th\>+\<a,(Sh\o)ch\t\>-\eps(h)\<a,c\>
=x_{a\t}\<h,(Sa\o)a\th\>}
where the last two terms cancel because $c$ is central. In other words,
the map $A\to H$ sending $a\mapsto (a\tens\id)\Delta c$ intertwines the
quantum adjoint action and the quantum coadjoint action; likewise for its
projection to $\ker\eps$, which is the map $a\mapsto x_a$. Using this
observation, we have the quantum Lie bracket
\[ [x_a,x_b]=x_{b\t}\<x_a,(Sb\o)b\th\>=x_{b\t}\<a,c\o\>\<c\t,(Sb\o)b\t\>
-x_b\<a,c\>\]
giving the result as stated. Likewise, the braiding in Proposition~2.4
comes out as
\align{\Psi(x_a\tens x_b)\equad &&=\<a,c\o\>\Ad_{c\t}(x_b)\tens c\th -
\<a,c\o\>\Ad_{c\t}(x_b)\tens 1\\
&&=x_{b\t}\tens c\t\<a(Sb\o)b\th,c\o\>-x_{b\t}\tens
1\<a(Sb\o)b\th,c\>=x_{b\t}\tens x_{a(Sb\o)b\th},}
again using the intertwining property for the map $a\mapsto x_a$.  \endproof

The centrally generated calculi are dual in a certain sense to inner type-I
calculi.
Thus, the quantum tangent space for the inner calculus in Proposition~2.5 can
be viewed as the kernel under differentiation for a suitable (right-handed)
calculus on $H$, taken  along the direction of $\alpha-\eps(\alpha)$. By
contrast, the quantum tangent space for a centrally generated calculus can be
viewed as the projection to $\ker\eps\subset H$ of the image of $c$ under
differentiation along all possible $a\in \ker\eps\subset A$, for a suitable
(left-handed) calculus on $H$. Also, for factorisable quantum groups (see the
next section) there is a correspondence between central elements and elements
invariant under a right handed $\Ad_R$ coaction, via the quantum Killing form,
see\cite{BrzMa:bic}. So centrally generated calculi and (a right-handed version
of) inner calculi are in correspondence in this case, although quite different
in character.

A centrally generated calculus typically has many quotients, i.e. its quantum
tangent space $L_c$ itself has further subrepresentations. For example,  we can
restrict the allowed $\{x_a\}$ to $a\in M_R$ whenever $M_R$ is a right ideal in
$\ker\eps\subset A$ stable under the right coaction $\Ad_R(a)\equiv a\t\tens
(Sa\o)a\th$. This is because the relation $ha=\<h,(Sa\o)a\th\>a\t$ holds for
the quantum double when acting on a central element. Moreover, comparing with
(\ref{quotcross}), we see that every non-trivial central element $c$ defines a
`mirror' operation from the moduli space of  right-handed bicovariant calculi
to  the moduli space of left-handed bicovariant calculi. It sends a calculus
defined by quotienting the universal one by  a right ideal $M_R$ according to a
(right-handed version of) (\ref{quotcross}) to the calculus with quantum
tangent space $\{x_a|\ a\in M_R\}\subseteq L_c$. The latter calculus is itself
a quotient of the universal one, namely by the (left) ideal
\eqn{mirror}{ \overline{ M_R}=\{a\in\ker\eps|\ \<ma,c\>=0\quad \forall m\in
M_R\}.}
This is our `mirror' operation at the level of the quotienting ideals. The same
operation with left-right interchanged takes us from left-handed to right
handed calculi, and there is a canonical inclusion
$M_R\subseteq\overline{\overline {M_R}}$. The calculus with quantum tangent
space $L_c$ in Proposition~2.6 is the mirror image of the zero differential
calculus, and vice-versa. In addition, the zero differential calculus is the
mirror image of the universal differential calculus.

Also clear from this point of view, if $L_1,L_2$ are subrepresentations of the
quantum double then $L_1\cap L_2$ is also. We denote its calculus by
$\Gamma_1\cdot\Gamma_2$; it is a quotient of both $\Gamma_1,\Gamma_2$. $L_1+
L_2$ is also a subrepresentation and we denote its calculus
$\Gamma_1*\Gamma_2$; it has $\Gamma_1,\Gamma_2$ as quotients. If $L_1\cap
L_2=\{0\}$ then the resulting calculus is  the obvious {\em direct product}
calculus. We say that a differential calculus is {\em  coirreducible}  if its
corresponding  quantum tangent space $L$ is irreducible as a quantum double
representation. This implies that the calculus has no proper
quotient calculus. Note that this should not be confused with irreducibility
for the calculus (no proper
subcalculus) which is automatically true for all calculi in the paper as a
consequence of axiom 5. in the definition of a bicovariant calculus.  Moreover,
in nice cases (where the quantum double is semisimple) one has only to
decompose
\eqn{kerPW}{ \ker\eps=L_1\oplus L_2\oplus \cdots }
into irreducibles in order to classify coirreducible calculi; each distinct
irreducible in the decomposition corresponds
to an isolated coirreducible calculus and each irreducible with multiplicity
typically corresponds to
a continuous family of calculi given by a parameter describing the embeddings
of the irreducible into its multiple copies  in $\ker\eps$. Moreover, we see in
this situation that
the universal differential calculus, which corresponds to $L=\ker\eps$, can be
built up as a direct product of coirreducible calculi.

\section{Calculi on finite groups and enveloping algebras}

In this section we apply the formulation of the classification problem in the
preceding section to the elementary cases of finite groups and enveloping
algebras.
The result in the case $A=\C(G)$ (the algebra of functions on a finite group)
is already known by other means, but recovered now from Proposition~2.2. But we
also give the dual case $A=\C G$ (the group algebra). It turns out to be more
similar to the quantum group case in the following section. The classification
for Lie groups and enveloping algebras remains open, but we make some remarks.

\begin{propos}  Let $A=\C(G)$ where $G$ is a finite group. It is known in this
case\cite{BMD:non} that
the coirreducible bicovariant differential calculi are in 1-1 correspondence
with the non-trivial conjugacy classes $\CC\subset G$. We recover this result
from the above approach as corresponding to
\[ L={\rm span}\{x_g\equiv g-e|g\in\CC\},\quad \del_{x_g} a=a(g\cdot(\
))-a,\quad \Psi^{-1}(a\tens x_g)=x_g\tens a(g\cdot(\ )) \]
\[ [x_g,x_h]=x_{ghg^{-1}}-x_h,\quad \Psi(x_g\tens x_h)=x_{ghg^{-1}}\tens  x_g\]
where $e$ is the group identity element.
\end{propos}
\proof Here $H=\C G$ and we classify all irreducible subspaces $L\subseteq
\ker\eps\subset \C G$ which are stable under the adjoint action and the action
of $\C(G)$ in Lemma~2.1. The algebra $A=\C(G)$ is commutative and elements of
the form $x_g=g-e$ are a basis of simultaneous eigenfunctions for its action on
$\ker\eps$, since $a\la x_g=a(g)g-a(e)e-a(g)e+a(e)e=a(g)x_g$ for any $g\in G$.
By choosing $a$ a Kronecker delta function we see that if  $L$ contains a
linear combination involving $x_g$ then it contains $x_g$ itself. Hence $L={\rm
span}\{x_g|g\in \CC\}$ for some subset $\CC\subset G$ not containing $e$. This
is the content of stability under the $A$ part of the quantum double action.
The content of stability under the $H$ part of the quantum double action (the
adjoint action of $G$ extended linearly)  is therefore that $\CC$ should be a
union of non-trivial conjugacy classes. The irreducible $L$ then correspond
precisely to the non-trival conjugacy classes. The corresponding braided
derivations are $\del_{x_g}a=\<x_g,a\o\>a\t=a(g(\ ))-a(e(\ ))$ and the `quantum
Lie bracket' is $[x_g,x_h]=\Ad_g(x_h)-x_h=x_{ghg^{-1}}-x_h$ as stated.
Likewise, we compute $\Psi$ from (\ref{Psixa}) and (\ref{Psixy}) in the form
stated. One also has $\ker\eps=\oplus_{\CC\ne\{e\}} L_{\CC}$, corresponding to
the decomposition of $G-\{e\}$
into non-trivial conjugacy classes, i.e. the universal calculus as a direct
product of the coirreducibles.
\endproof

These calculi are all `non-classical' in the sense that the braiding needed for
the derivation property is non-trivial (when $G$ is non-Abelian). They are in
fact a variant of the familiar $q$-derivative, with $q$ being replaced by a
group element taken from the conjugacy class.
The non-classical nature also appears as non-commutativity of the calculus in
the sense $a\extd b\ne (\extd b) a$ for some $a,b$.
The calculus is  inner type-I with $\alpha$ the characteristic function of
$\CC\cup\{e\}$, and inner type-II with $\alpha$ the characteristic
function of $\CC$. We can also apply our formalism  to $A=\C G$. If
$G$ is Abelian we have $\C G=\C(\hat G)$ and return to the preceding example
applied to the dual group. But when $G$ is non-Abelian, the algebra $A$ is
non-commutative and we are really doing `non-commutative geometry'.

\begin{propos} Let $A=\C G$ where $G$ is a finite group. The coirreducible
bicovariant differential calculi
are in 1-1 correspondence with pairs $V,\lambda$, where $V$ is a non-trivial
irreducible (right)   representation of $G$ and
$\lambda\in P(V^*)$. The corresponding calculus has dimension $\dim V$ and
\[ L={\rm span}\{x_v\equiv \<v\ra(\ ), \lambda\>-\<v,\lambda\>1\ |\ v\in
V\},\quad \del_{x_v}(g)=(\<v\ra g,\lambda\>-\<v,\lambda\>)g,\quad
\Psi^{-1}(x_v\tens g)=g\tens x_{v\ra g}\]
and trivial `quantum Lie bracket'.
\end{propos}
\proof Here $H=\C(G)$ is commutative. Hence the adjoint action in Lemma~2.1 is
trivial (as is the
bracket $[\ ,\ ]$ and its associated braiding). We therefore need only to
classify irreducible subspaces $L\subseteq\ker\eps\subset \C(G)$ under the
action
of $\C G^{\rm op}$. This action is $h\la x=x(h(\ ))-x(h)1$ for all
$x\in\ker\eps$, which is the standard projection $\Pi$ to $\ker\eps$ of the
right regular representation of $G$
on $\C(G)$ by multiplication from the left in the group. The Peter-Weyl
decomposition $\C(G)\isom \C\oplus_{V\ne \C} V\tens V^*$ projected via the
projection $\Pi$ is an isomorphism $\oplus_{V\ne \C}V\tens V^*\isom \ker\eps$,
giving the decomposition of this into irreducibles. In the Peter-Weyl
decomposition, the element $v\tens \lambda$
maps to the function  $\<v\ra(\ ),\lambda\>\in \C(G)$, giving the form of $L$
shown. We need $\lambda\ne 0$ and we identify all $\lambda$ which are related
by a phase since these give the same $L$, i.e. the continuous parameter is
$\lambda\in P(V^*)=\C P^{\dim V-1}$. The braided-derivation is
$\del_{x_v}g=\<x_v,g\>g$ on group-like elements of $\C G$, which gives the form
shown. The group-like
elements are simultaneous eigenfunctions for  all the braided-derivations. The
braiding is easily computed
as $\Psi^{-1}(g\tens x_v)=  g\la x_v\tens g = \<v\ra g(\ ),\lambda\>\tens g-
\<v\ra g,\lambda\>\tens g=  x_{v\ra g}\tens g$.
\endproof

Note that a basis of $V^*$ specifies $\dim V$ isomorphic copies of $V$ in the
Peter-Weyl decomposition. However,
we need here not only the multiplicities but the actual corresponding subspaces
$L$. We obtain a
subspace isomorphic to $V$ for every non-trivial linear combination (modulo an
overall scale) of the basis elements, i.e. a continuous family
of calculi parametrized by the projective space $P(V^*)$ for each irreducible
representation $V$. Also, since irreducible representations of $G$ correspond
to characters, one can recast this result in terms of these.
For a given character $\chi$ we identify $V_\chi^*$ as the quotient of $\C G$
where $[\lambda]=[\lambda']$ if $\chi(g\lambda)=\chi(g\lambda')$ for all $g$.
Then coirreducible calculi are in 1-1 correspondence with
pairs $\chi,[\lambda]$ according to
\[ L={\rm span}\{x_g\equiv \chi(g(\ )\lambda)-\chi(g\lambda)1|g\in G\},\quad
\del_{x_g}h=(\chi(gh\lambda)-\chi(g\lambda))h,\quad \Psi^{-1}(h\tens x_g )=
x_{gh}\tens h\]
where $g,h\in G$. Here $V_\chi$ is the vector space spanned by $\chi((\
)\lambda g)$ as $g$ runs over $G$ and  is
an irreducible (right) representation of $G$ acting by left multiplication in
the argument of $\chi$. From this, it is clear that these calculi on $\C G$
are all centrally generated by $c=\chi(g(\ )\lambda)$.

Finally, we consider the differential calculi on a classical Lie group
co-ordinate ring $A=\C(G)$.  Here $\C(G)$ denotes an algebraic model of the
functions on $G$ constructed as a Hopf algebra non-degenerately paired to
the enveloping algebra $U(\cg)$.

\begin{propos} Let $\cg$ be a Lie algebra. For each natural number $n$ there
is a bicovariant differential calculus with
$L=\cg+\cg\cg+\cdots\cg^n$, the subspace of degree $\le n$ and $\ge 1$. For
example, for $n=2$:
\[ L={\rm span}\{ \xi, \eta\zeta|\ \xi,\eta,\zeta\in \cg\},\quad
\del_\xi=-\tilde\xi,\quad \del_{\eta\zeta}=\tilde\zeta\tilde\eta\]
\[\Psi^{-1}(a\tens \xi )= \xi\tens a,\quad \Psi^{-1}(a\tens \eta\zeta )=
\eta\zeta\tens a-\zeta\tens \tilde\eta(a)  - \eta\tens\tilde\zeta(a) \]
\[ [\xi,x]=\xi x-x\xi,\quad [\eta\zeta,x]=\eta\zeta x -\eta x \zeta - \zeta x
\eta + x\zeta\eta\]
\[\Psi(\xi\tens x)=x\tens\xi,\quad \Psi(\xi\eta\tens
x)=[\xi,x]\tens\eta+[\eta,x]\tens \xi+x\tens\xi\eta\]
for all $x\in L$. Here $\tilde\xi$ is the right-invariant vector-field
associated to $\xi\in\cg$.
\end{propos}
\proof  Note that the degree of
a given element in $U(\cg)$ is not well-defined but the subspace spanned by
products of up to $n$ elements is. We show only that such a subspace $L^{(n)}$
forms a quantum double subrepresentation. To see that it is closed under the
adjoint action of $U(\cg)$ it suffices to see that it is closed under the
action of each $\xi\in\cg$. This action is by commutator in $U(\cg)$. Hence
assuming the result for $L^{(n-1)}$ and the Leibniz rule for commutators, we
obtain the result for $L^{(n)}$ by induction. The other part of the quantum
double action (that of $A=\C(G)$) is given by evaluation against the left
coaction $\beta=\Delta -\id\tens 1$.
Then $\beta(\xi x)=\Delta(\xi x)- \xi x\tens 1=(\xi\tens 1)\Delta x+
(1\tens\xi)\Delta x- \xi x\tens 1=(\xi\tens 1)\beta(x)+(1\tens
\xi)\beta(x)+ x\tens \xi\in U(\cg)\tens L^{(n)}$ as $\beta(x)\in U(\cg)\tens
L^{(n-1)}$ and $n\ge 1$. Here $x\in L^{(n-1)}$ and we proceed by induction.
The explicit
computations for $L^{(2)}$ are immediate from the form of the coproduct on
$\eta\zeta$ in the formulae above.
Here, $\del_\xi(a)=\<\xi,a\o\>a\t={d\over dt}|_0 a(e^{t\xi}(\ ))=-\tilde\xi(a)$
for $a\in\C(G)$ (this is given explicitly by the matrix representation of $\cg$
used in defining the pairing between $U(\cg)$ and $\C(G)$, i.e. it is actually
algebraic.) \endproof

We see that it is possible to view higher order differential operators as if
they are `first order vector fields' -- but braided. A second order operator,
for example, is clearly not a derivation in the usual sense but it is a
braided-derivation for suitable $\Psi$. For example, one could compute its
`flow' as a corresponding braided-exponential. This opens up the possibility of
a `geometrical' picture for the evolution of quantum systems generated by
second or higher order Hamiltonians, to be given in detail elsewhere.

On the other hand,  we do not attempt to classify all bicovariant calculi here.
This would appear to be an interesting problem in
the classical theory of enveloping algebras: find all subspaces $L$ which are
stable under the adjoint action and under the left
coaction $\beta=\Delta-1\tens\id$. Moreover, the $L^{(n)}$ are of course not
coirreducible. Instead, we have a filtration
\eqn{filt}{ \cg=L^{(1)}\subset L^{(2)}\subset L^{(3)}\cdots\subset
L^{(\infty)}=\ker\eps,}
where $\cg=L^{(1)}$ corresponds to the classical differential calculus on
$\C(G)$. At the level of bicovariant calculi we have a sequence of quotients of
the universal one (of all finite degree invariant differential operators)
eventually quotienting down to the standard one.

There are certainly bicovariant calculi other than the $L^{(n)}$. For example,
if $\cg\tens\cg$ has an $\Ad$-invariant element $t=t_i\tens t^i$ (e.g. if $\cg$
is semisimple) then $L=\cg\oplus\C$ spanned by $\cg$ and the  central element
$c=t_it^i$
corresponds to a bicovariant differential calculus in between those
corresponding to $L^{(1)}$ and $L^{(2)}$. In the semisimple case, $\cg$ can be
viewed as being centrally generated according to Proposition~2.6 by the
quadratic Casimir, and $\cg\oplus \C$ is its canonical extension.
(Equivalently, the mirror operation (\ref{mirror}) in this case turns the zero
differential calculus into the  classical one, and vice-versa.) The elements of
$\cg$ act as ordinary vector fields, while $c$ acts as a second order operator
viewed as a braided-vector field. The quantum Lie bracket
restricted to $\cg$ is its usual Lie bracket. The other cases and the braiding
are
\ceqn{quadLc}{ [\xi,c]=0,\quad [c,c]=0,\quad [c,\xi]=[t_i,[t^i,\xi]],\quad
 \Psi(\xi\tens\eta)=\eta\tens\xi\\
 \Psi(\xi\tens c)=c\tens\xi,\quad \Psi(c\tens\xi)=\xi\tens c+[t_i,\xi]\tens
t^i-t^i\tens[t_i,\xi],\quad \Psi(c\tens c)=c\tens c.}
If $\cg$ under the adjoint action is isotypical (as for $sl_2$) then $[c,\xi]$
here is fixed multiple of $\xi$. The simplest case $L=sl_2\oplus \C$
corresponds to the non-standard 4-dimensional differential calculus on $SU(2)$
which has been studied in \cite{BMD:non} as the $q\to 1$ limit of the  known
4-dimensional calculus on the quantum group $SU_q(2)$ in \cite{Wor:dif}.
Similarly, $L^{(n-1)}\oplus \C$  corresponds to a natural calculus in between
the calculi corresponding to $L^{(n-1)}$ and $L^{(n)}$, whenever we have a
degree $n$ central element. Intermediate calculi are generally what arise when
we take the limit of quantum group differential calculi (these will be
classified in the next section), i.e. this is a general feature. Put another
way, we will see
from the classification in the next section that the standard
$\dim\cg$-dimensional calculus on a simple Lie group $G$ violates the
`principle of $q$-deformisability'; only certain extensions of ordinary vector
fields on a Lie group by higher order
vector fields can deform to calculi on $G_q$.

\section{Calculi on factorisable quantum groups}

In this section we present our main result, which is a classification of the
bicovariant calculi for a factorisable
semisimple quantum group. We then discuss the application of the result to the
standard quantum groups $G_q$, which these are essentially factorisable.

We recall that a `strict quantum group' or quasitriangular Hopf algebra is
factorisable\cite{ResSem:mat} if   $\CR_{21}\CR$ viewed as a map $\CQ:A\to H$
by
 $\CQ(a)=(a\tens\id)(\CR_{21}\CR)$ is an isomorphism. This is the strongest
form;
one may also
demand separately that the map is injective or surjective. We also consider, by
definition, that a quantum group is semisimple if there is a Peter-Weyl
decomposition
\eqn{PW}{   \oplus_{V} V^*\tens V\isom A}
provided by the matrix elements of the inequivalent finite-dimensional
irreducible representations $V$ of $H$.
This is broadly equivalent to other notions of semisimplicity, and is at any
rate the condition that we suppose in this section. If $V$ is such a
representation, with basis $\{e_i\}$ and dual basis $\{f^i\}$, we define the
matrix elements $\rho^i{}_j\in A$ by $h\la e_i=e_j \rho^i{}_j(h)$
and the above map by $f^i\tens e_j\mapsto \rho^i{}_j$.

\begin{lemma} Let $H$ be a factorisable quantum group with dual $A$. The map
$\CQ$ identifies $\ker\eps\subset A$ and $\ker\eps\subset H$. Under this
identification, the action of the quantum double in Lemma~2.1 becomes the
action on $\ker\eps\subset A$ given by
\[ h\la a=a\t \<h,(Sa\o)a\th\>,\quad b\la a= \<b,\CR'\uo
\CR\ut\>\<a\o,\CR'\ut\>\<a\th,\CR\uo\>  a\t -\<b, \CQ(a)\>  1\]
for all $h\in H$ and $b\in A$. Here, $\CR'\equiv\CR'\uo\tens\CR'\ut$ denotes a
second copy of  the quasitriangular structure.
$\CR$.
\end{lemma}
\proof It is immediate from the counity property of the quasitriangular
structure $\CR$ that $\CQ(1)=1$. Hence $\CQ(\ker\eps)=\ker\eps\subset H$.
Moreover,
we know from $\Ad$-invariance of the quantum Killing form that $\Ad_h \circ
\CQ(a)=\CQ(\Ad^*_h a)$ where $\Ad_h^*$ is the left quantum coadjoint
action as
stated for $h\la a$ in the lemma, and $\Ad$ is the quantum adjoint action used
for $x$ in Lemma~2.1. For a proof see \cite{Ma:skl} or the text\cite{Ma:book}.
The new part concerns the other action:
\align{b\la \CQ(a)\equad &&=\<b,\CQ(a)\o\>\CQ(a)\t - \<b,\CQ(a)\>1\\
&&=\<a,\Rut1\Ruo2\>\<b,\Ruo1\o\Rut2\o\>\Ruo1\t\Rut2\t-\<b,\CQ(a)\>\CQ(1)\\
&&=\<a,\Rut1\Rut3\Ruo4\Ruo2\>\<b,\Ruo1\Rut2\>\Ruo3\Rut4-\<b,\CQ(a)\>\CQ(1)\\
&&=\<a\o,\Rut1\>\<a\th,\Ruo2\>\<b,\Ruo1\Rut2\>\CQ(a\t)-\<b,\CQ(a)\>\CQ(1)
=\CQ(b\la a)}
for all $a\in\ker\eps\subset A$ and $b\in A$, where
$\Ruo1\tens\Rut1,\cdots,\Ruo4\tens\Rut4$ are four copies of $\CR$. The first
equality is the action of $A$ in Lemma~2.1. The second puts in the formula for
$\CQ$. The third is the coproduct property of the quasitriangular structure
and finally we recognise the required result in terms of the action $b\la a$
stated. Hence $\CQ$ intertwines the stated action of the quantum double with
the action in Lemma~2.1. Note that this computation also works at the level of
a
coaction of $H$ rather than an action by $b\in A$ (i.e. the action of the
quantum double remains $A$-regular). \endproof

So the possible quantum tangent spaces $L$ are in 1-1 correspondence with
subrepresentations
of $\ker\eps\subset A$ under this action of the quantum double. This action
looks more complicated than before. However, there is a well-known isomorphism
in the factorisable case of the quantum double with $H\codcross H$. The latter
is $H\tens H$ as an algebra and has a coalgebra which is a twisting of the
tensor product one. The map $\theta$  to $H\codcross H$  is\cite{ResSem:mat}
\ceqn{facisom}{ \theta(h\tens a)=h\o \CR\umt\tens h\t\CR\uo\<\CR\umo\CR\ut,a\>}
The full details of the isomorphism and an explicit formula for $\theta^{-1}$
are in the author's text\cite{Ma:book}.

\begin{propos} The action in Lemma~4.1 of the quantum double, in the form
$H\codcross H$ acting on $\ker\eps\subset A$, takes the form
\[ (h\tens 1)\la a=\<Sh,a\o\>a\t - 1\<Sh,a\>,\quad (1\tens g)\la
a=a\o\<g,a\t\>-1\<g,a\>\]
for all $h,g\in H$ and $a\in\ker\eps\subset A$.
\end{propos}
\proof To find the action of $H\codcross H$ we need the explicit inversion
formula for $\theta$ in \cite{Ma:book}. Then $(h\tens 1)\la
a=\theta^{-1}(h\tens 1)\la a$ etc. can be computed, and one obtains the result
stated in the proposition. Once these actions have been obtained, however, it
is   enough (and rather easier) to verify that pull back along $\theta$ indeed
recovers the action of $H\dcross A^{\rm op}$ in Lemma~4.1. Thus,
\align{\theta(h\tens 1)\la a\equad &&=(h\o\tens h\t)\la a=(h\o\tens 1)\la
a\o\<h\t,a\t\>-(h\o\tens 1)\la 1\<h\t,a\>\\
&&=\<Sh\o,a\o\>a\t\<h\t,a\th\>-\<Sh\o,a\o\>1\<h\t,a\t\>=\<h,(Sa\o)a\th\>a\t\\
\theta(1\tens b)\la a\equad &&=(\CR\umt\tens \CR\uo)\la a\<\CR\umo\CR\ut,b\>\\
&&=(\CR\umt\tens 1)\la a\o \<\CR\uo,a\t\>\<\CR\umo \CR\ut,b\>-(\CR\umt\tens
1)\la 1\<\CR\uo,a\>\<\CR\umo\CR\ut,b\>\\
&&=\<S\CR\umt,a\o\>a\t\<\CR\uo,a\th\>\<\CR\umo\CR\ut,b\>
-\<S\CR\umt,a\o\>1\<\CR\uo,a\t\>\<\CR\umo\CR\ut,b\>\\
&&=a\t\<\CR'\ut,a\o\>\<\CR\uo,a\th\>\<\CR'\uo\CR\ut,b\>-\<\CQ(a),b\>1}
as required. We used the form of $\theta$, the actions as stated in the
proposition and, in the last line, the antipode property
$(S\tens\id)\CR^{-1}=\CR$ of a quasitriangular structure.
Our notation is $\CR\umo\tens \CR\umt=\CR^{-1}$.
\endproof

So, quantum tangent spaces $L$ are in correspondence with subrepresentations of
$\ker\eps$ under this action of $H\tens H$. We can now obtain our main result.

\begin{theorem} Let $H$ be a factorisable quantum group with dual $A$, and
suppose that the Peter-Weyl decomposition (\ref{PW}) holds. Then the finite
dimensional bicovariant coirreducible calculi on $A$ are in 1-1 correspondence
with the non-trivial finite-dimensional irredicible representations $V$ of $H$.
The corresponding calculus has dimension $(\dim V)^2$ and
\[ L={\rm span}\{x^i{}_j\equiv \CQ(\rho^i{}_j-1\delta^i{}_j)|\ i,j=1,\cdots,
\dim V\}\]
\[ \del_{x^i{}_j}(a)=\CQ(\rho^i{}_j\tens a\o)a\t-\delta^i{}_j a,\quad
\Psi^{-1}(a\tens x^i{}_j)=x^a{}_b\tens a\th
\CR(a\o\tens\rho^i{}_a)\CR(\rho^b{}_j\tens a\t)\]
\[ [x^i{}_j,x^k{}_l]=x^a{}_b\CQ(\rho^i{}_j\tens
(S\rho^k{}_a)\rho^b{}_l)-x^k{}_l\delta^i{}_j\] \[\Psi(x^i{}_j\tens
x^k{}_l)=x^m{}_n\tens
x^a{}_b\CR((S\rho^c{}_m)\rho^n{}_d\tens\rho^i{}_a)\CR(\rho^b{}_j
\tens(S\rho^k{}_c)\rho^d{}_l)\]
where we also regard the quantum Killing form and quasitriangular structure
as functionals $\CQ,\CR:A\tens A\to \C$
\end{theorem}
\proof We first separate off the trivial representation in (\ref{PW}), so
$A\isom \C\oplus(\oplus_{V\ne\C} V^*\tens V)$ where the sum is over non-trivial
$V$. The projection $\Pi(a)=a-1\eps(a)$ from $A\to \ker\eps$ establishes an
isomorphism
\eqn{kerV*V}{ \ker\eps\isom \oplus_{V\ne\C}V^*\tens V.}
This is because $\Pi$ and the projection to $\oplus_{V\ne\C}V^*\tens V$ have
the same kernel, namely the span of the identity element in $A$. By
Proposition~4.2, we therefore have an isomorphism of $H\tens H$ modules, where
the second $H$ acts on $V$ as in the
Peter-Weyl decomposition (the given irreducible representation $V$) and the
first copy of $H$ acts on $V^*$ by the conjugate representation $h\la
f=f(Sh\la(\ ))$ for $f\in V^*$. Next, as $H\tens H$ modules, these $V^*\tens V$
are distinct and irreducible. Hence they are precisely the choices for
irreducible subrepresentations of $\ker\eps\subset A$.

The explicit formula for the braided-derivations and their requisite braiding
are easily computed from the formulae in Proposition~2.3.  From the proof of
Lemma~4.1 we have
\align{(\Delta-\id\tens 1)x^i{}_j\equad
&&=\CR\uo\CR'\ut\<(\rho^i{}_j-\delta^i{}_j)\o,\CR\ut\>\<(\rho^i{}_j
-\delta^i{}_j)\th,\CR'\uo\>\tens \CQ((\rho^i{}_j-\delta^i{}_j)\t)
-\CQ(\rho^i{}_j-\delta^i{}_j)\tens 1\\
&&=\CR\uo\CR'\ut\tens\<\rho^i{}_a,\CR\ut\>\<\rho^b{}_j,\CR'\uo\>
\CQ(\rho^a{}_b)-\CQ(\rho^i{}_j)\tens 1\\
&&=\CR\uo\CR'\ut\tens\<\rho^i{}_a,\CR\ut\>\<\rho^b{}_j,\CR'\uo\>x^a{}_b}
Evaluation against this is the action of $A$ in Lemma~4.1, which is the action
needed to compute the braiding. Thus, $\Psi^{-1}(a\tens x^i{}_j)=a\t\tens
\<a\o,\CR\uo\CR'\ut\>\<\rho^i{}_a,\CR\ut\>\<\rho^b{}_j,\CR'\uo\>x^a{}_b$, which
can be written in the form shown where  $\CR$ is regarded as a functional on
$A\tens A$. The quantum Lie bracket and its braiding from Proposition~2.4 are
also easily computed and follow the same lines as in
\cite{Ma:skl}\cite{Ma:exa}, except that we are not tied to any particular
representation $V$ or any fixed R-matrix; we include the proofs only for
completeness in our present conventions. Thus, by $\Ad$-invariance of $\CQ$ we
have
\align{[x^i{}_j,x^k{}_l]\equad
&&=\CQ[\CQ(\rho^i{}_j-\delta^i{}_j)\la(\rho^k{}_l-\delta^k{}_l)]\\
&&=\CQ(\rho^a{}_b)\<\CQ(\rho^i{}_j),(S\rho^k{}_a)\rho^b{}_l\>-\delta^i{}_j
\CQ(\rho^k{}_l)
=x^a{}_b\<\CQ(\rho^i{}_j),(S\rho^k{}_a)\rho^b{}_l\>-\delta^i{}_jx^k{}_l,}
which we write in the form stated where $\CQ=\CR_{21}\CR$ is regarded as a
functional on $A\tens A$. Here $\la$ is the quantum coadjoint action of $H$
in Lemma~4.1. Finally, using the above result for $\Delta x^i{}_j$ and
$\Ad$-invariance of $\CQ$, we have
\align{\Psi(x^i{}_j\tens x^k{}_l)\equad &&=[x^i{}_j\o,x^k{}_l]\tens
x^i{}_j\t-[x^i{}_j,x^k{}_l]\tens 1=\CQ[\CR\uo\CR'\ut\la(\rho^k{}_l
-\delta^k{}_l)]
\tens\<\rho^i{}_a,\CR\ut\>\<\rho^b{}_j,\CR'\uo\>x^a{}_b\\
&&=\CQ(\rho^c{}_d)\tens x^a{}_b\<\CR\uo\CR'\ut,(S\rho^k{}_c)\rho^d{}_l\>
\<\rho^i{}_a,\CR\ut\>\<\rho^b{}_j,\CR'\uo\>x^a{}_b-\delta^k{}_l
\tens x^i{}_j\\
&&=x^c{}_d\tens x^a{}_b\<\CR\uo\CR'\ut,(S\rho^k{}_c)\rho^d{}_l\>
\<\rho^i{}_a,\CR\ut\>\<\rho^b{}_j,\CR'\uo\>x^a{}_b}
which we write in the form stated. Note that both of the expressions
$\CQ(\rho^i{}_j\tens (S\rho^k{}_a)\rho^b{}_l)$ and
$\CR((S\rho^c{}_m)\rho^n{}_d\tens\rho^i{}_a)\CR(\rho^b{}_j\tens(S\rho^k{}_c)
\rho^d{}_l)$ can be expanded out as four-fold products of the matrices
$R=(\rho\tens\rho)\CR$, its inverse and $\tilde R
=(\rho\tens\rho\circ S)\CR$. This step and the resulting R-matrix formulae
are identical in form to the computation of the quantum Lie algebra
`structure constants' in \cite{Ma:skl} and the computation of the quadratic
relations of the braided matrices in \cite{Ma:exa} (the matrix denoted
$\Psi'$ there), respectively. Hence we omit the  proofs and note only
that, after rearranging the R-matrices, one has the same form as
for a quantum or braided-Lie algebra of matrix type, namely
\eqn{qlieR}{ R_{21}[x_1,Rx_2]=x_2Q-Qx_2,\quad R_{21}\Psi(x_1\tens
Rx_2)=x_2R_{21}\tens x_1R}
where the numerical suffices denote positions in a matrix tensor product and
$Q=R_{21}R$. The relation between (\ref{qlieR}), braided matrices
$\vecu=x+\id$ and the quantum double is
explained further in \cite{Ma:sol} (where the quantum double braiding $\Psi$ is
denoted $\haj{{\bf R}}$). On the other hand, now (\ref{qlieR})
applies to any irreducible representation $V$ of $H$ and not some fundamental
basic representation, which need not exist. \endproof

Let us note that if $\CR$ is a quasitriangular structure in a quantum group
then so is $\CR^{-1}_{21}$. Thus all results involving a quasitriangular Hopf
algebra have a `conjugate' one in which this conjugate $\CR_{21}^{-1}$ is used
instead of $\CR$. This conjugation is also intimately tied to the $*$-operation
or complex conjugation in many systems\cite{Ma:qsta}. In the above theorem, we
see that for every $V$ we have equally well the conjugate
\eqn{barL}{  \bar L={\rm span}\{\bar x^i{}_j\equiv \bar{
\CQ}(\rho^i{}_j-1\delta^i{}_j)|\ i,j=1,\cdots,\dim V\}}
where $\bar {\CQ}(a)=(a\tens\id)(\CR^{-1}\CR^{-1}_{21})$. Here $\bar L$ is
isomorphic   to $L$ but the isomorphism (which is $\bar {\CQ}\circ \CQ^{-1}$
restricted to $L$) is  non-trivial. This fits also with the general point of
view of quasi-$*$ structures on inhomogeneous quantum groups\cite{Ma:qsta}
where the tensor product of unitaries is unitary only up to a non-trivial
isomorphism.

These results can be applied formally to the standard quantum groups
$H=U_q(\cg)$ with dual $A=G_q$ associated to complex semisimple Lie algebras,
provided we work over formal power-series $\C[[\hbar]]$ and introduce suitable
logarithms for some of the $G_q$ generators, etc. Or, if we want to work
algebraically over $\C$ (with generic $q$), we need to localise and introduce
roots of some of the generators of $G_q$ and use the algebraic form of
$U_q(\cg)$ where $q^{H\over 2}$ etc. is regarded as a single generator. This is
clear from the standard cases such $SU_q(2)$: In standard notations the value
of $\CQ$ on the generators is
\eqn{su2}{ \CQ\pmatrix{a&b\cr c&d}=\pmatrix{q^H&q^{-\h}(q-q^{-1})
q^{H\over 2}X_-\cr
q^{-\h}(q-q^{-1})X_+q^{H\over 2}& q^{-1}C-q^{-2}q^H}}
where $C=q^{H-1}+q^{-H+1}+(q-q^{-1})^2X_+X_-$ is the $q$-quadratic Casimir.
According to \cite{ResSem:mat}, the standard
quantum groups are all factorisable modulo such formal extensions. Likewise,
the Peter-Weyl decomposition (\ref{PW}) holds formally for the standard
semisimple $\cg$. This is because the category of finite-dimensional
representations in the classical and quantum cases are generically equivalent,
and the assumption holds in some form for the classical case.  Note also
that the
entries in (\ref{su2}) projected to $\ker\eps$ span a 4-dimensional $L$
associated to the spin $1/2$ representation, and has the structure in
Theorem~4.3 without any powerseries.  Indeed, the
quantum double of $U_q(su_2)$ is known to be a $q$-deformation of the
Lorentz group and hence the lowest possible generic representation
is the 4-dimensional one on $q$-Minkowski space. In this simplest case,
$\bar L$ is the same subspace $L$. The latter also coincides with $L_C$ from
Proposition~2.6 with $C$ the $q$-quadratic Casimir above, and is a subspace of
$L_{\alpha,1}$ from Proposition~2.5 with
$\alpha=(qa+q^{-1}d)/q^{-2}(q^3-1)(q-1)$ the normalised $q$-trace.

Therefore we should understand Theorem~4.3 not as a complete algebraic
classification for a given version of each given $G_q$ (this is a much harder
problem and has been recently addressed in some cases\cite{SchSch:bic}), but as
a classification of those calculi which are `generic' in the sense that they
extend to the various localisations and square-roots of the generators etc.
needed for exact factorisability. In other words, there are  natural calculi,
corresponding to  $L$ (or $\bar L$) for each finite-dimensional irreducible
representation $V$, and these are the only ones modulo `pathological'
possibilities for particular $q$ for particular versions of particular $G_q$.
\note{Moreover, we are now able to understand the well-known phenomenon that
some of the calculi in the case-by-case classification, which is known for some
quantum groups, come in pairs, such as the $4D_\pm$ calculi for $SU_q(2)$ in
\cite{Wor:dif} and similarly for $SU_q(n)$ in \cite{SchSch:bic}. We see that
this is (a) a general feature arising from braid-crossing-reversal (which leads
to both $L,\bar L$) and (b) not visible at the level of formal power-series; at
this level $L\isom \bar L$ but the isomorphism $\bar {\CQ}\circ \CQ^{-1}$
(restricted
to $L$) does not have a polynomial algebraic description in terms of the
generators and product of $U_q(\cg)$, and therefore might be hard to find by
strictly algebraic methods.  }

For the $A,B,C,D$ series we have a natural `fundamental representations' $V$
and in this case it should be clear that the calculus corresponding to $L$ is
the one found by Jurco\cite{Jur:dif} by other means. We therefore have a new
construction for this and the result that its slight generalisation to other
irredcuble representations exhausts all the generic first order bicovariant
differential calculi on the standard semisimple quantisations.

\section{Concluding remarks}

We conclude with some remarks about further work. Firstly, the bicovariant
calculi studied here are `first order'. They play the role of 1-forms. It
remains to construct and classify all possible higher order calculi or
`exterior algebras'. One canonical construction is to take the tensor algebra
on the first order calculus $\Gamma$ and quotient with the aid of a
`skew-symmetrizer' built from the quantum double braiding $\Psi$,
see\cite{Wor:dif}. In this case the exterior algebra is a super-Hopf
algebra\cite{Brz:rem}. On the other hand, even when $\Gamma$ is the classical
calculus, the canonical exterior algebra is not the classical one. One must
quotient it further. The classification of exterior algebras therefore
remains open even after we have classified the first order calculi.

Secondly, all of the results in Section~2 about first order calculi on quantum
groups have an analogue for braided groups. Braided groups are needed to
include $q$-deformations $\R_q^n$ and $\R_q^{1,3}$ etc., with their additive
(braided) coproduct. The classification of differential calculi on such objects
would therefore seem to be the starting point for some form of $q$-geometry
based on $\R^n$. Our result in this direction is a negative but rather
unexpected one: generically there is only one coirreducible braided-bicovariant
differential calculus on $\R_q^{1,3}$ (say), and it is infinite-dimensional.
Its braided tangent space $L$ consists (in a suitable completion) of a
$q$-deformation of the space of solutions of the massless Klein-Gordon equation
projected to the functions vanishing at the origin.
Briefly, (details will be presented elsewhere) the sketch is as follows. Let
$B$ be a braided group in a braided category generated by `background quantum
group' $H$ as its category of modules. We define a braided-bicovariant calculus
$\Gamma$ in the obvious way and proceed in a similar manner to Section~2. The
role of the quantum double is now played by the author's `double-bosonisation'
$B^*\lbiprod H\rbiprod B$ quantum group\cite{Ma:dbos}. This acts on
$\ker\eps\subset B$ and the possible braided tangent spaces $L$ are in 1-1
correspondence with sub-representations of $\ker\eps$. When $B=\R_q^{1,3}$ it
is known from \cite{Ma:geo} that the double-bosonisation is the $q$-conformal
group and the action on $B$ is a q-deformation of its action on
$\R^n$. Classically, however, this representation has one irreducible
subrepresentation, which is the space of solutions of the massless Klein-Gordon
equation.

The braided version of the theory may also help to solve the above-mentioned
problem of exterior algebras on quantum groups, at least in the case of strict
(quasitriangular) quantum groups.  This is because the braided groups
corresponding under transmutation to strict quantum groups are always
braided-commutative in a certain sense\cite{Ma:bg}, i.e. closer to the
classical situation. Using this braided-commutativity one may reasonably expect
a natural exterior algebra q-deforming the classical one.  Such a result would
be the `skew' analogue of the situation in Section~2, where we explained that
the braided version of the `quantum Lie bracket' is better behaved for
constructing a braided enveloping algebra. This remains a direction for further
work.

\baselineskip 19pt


\end{document}